\documentclass[]{article}
\usepackage{latexsym}
\usepackage{amsmath,amsfonts,amssymb,theorem}
\newcommand{\R}{\mathbb R}
\newcommand{\Z}{\mathbb Z}

\newcommand{\C}{\mathbb C}
\newcommand{\Q}{\mathbb Q}

\newcommand{\kitem}{\begin{itemize}\vspace{-2ex}}
\newcommand{\kenditem}{\vspace{-1ex}\end{itemize}}
\newcommand{\CA}{{\cal A}}

\newcommand{\supp}{\mbox{\rm supp}}

\newcommand{\surj}{\rightarrow\hspace{-0.8em}\rightarrow}

\newcommand{\ko}{\overline}

\newcommand{\ku}{\underline}

\newcommand{\factor}[6]{{\raisebox{#3ex}{$#1$}\hspace{-#5em}\big/\hspace{-#6em}\raisebox{-#4ex}{$#2$}}}

\newcommand{\fac}[2]{\factor{#1}{#2}{0.4}{0.5}{0.3}{0.1}}
\newcommand{\pZ}{\Z_{\geq 0}}

\newcommand{\pR}{\R_{\geq 0}}

\newcommand{\kP}{P}
\newcommand{\height}{\mbox{\rm ht}}
\newcommand{\kk}{{\C}}

\newcommand{\qed}{{
\unskip\nobreak\hfil\penalty50\hskip0.1em\hbox{}\nobreak\hfil$\Box$
\parfillskip=0pt \finalhyphendemerits=0 \par\medbreak}}

\reversemarginpar
 \newcommand{\kkk}[1]{}

\newcounter{Abschnitt}[section]
\newcommand{\neu}[1]{\protect\refstepcounter{Abschnitt}\protect
   \label{#1}\vspace{1ex}
   {\bf (\protect\arabic{section}.\protect\arabic{Abschnitt})}
                     $\qquad$\kkk{#1}}
\newcommand{\zitat}[2]{{\rm(}\protect\ref{#1}.\protect\ref{#1-#2}{\rm)}}

\begin{document}
\title{\bf\huge The chain property for the associated primes of $\CA$-graded
ideals}
\author{Klaus Altmann\footnote{supported by a DFG Heisenberg grant} }
\date{}
\maketitle
\begin{abstract}
We investigate how the chain property for the associated primes of monomial
degenerations of toric (or lattice) ideals can be  
generalized to arbitrary $\CA$-graded ideals.
The generalization works in dimension $d=2$, but it fails for $d\geq 3$.
\end{abstract}

%
%
\section{Introduction}\label{intro}


\neu{intro-agraded}
Challenged by the question of Arnold for the ideals with the easiest
Hilbert function, Sturmfels has invented in \cite{agraded} and \S 10 of
\cite{polytopes} the notion of $\CA$-graded ideals. For a given linear
map $\CA:\Z^n\to\Z^d$ with $(\ker\CA)\cap\pZ^n=0$ 
an ideal $I\subseteq \kk[x_1,\dots,x_n]$ is called $\CA$-graded
if it is $\Z^d$-homogeneous via $\CA$ and, moreover, if it has the Hilbert 
function
\[
\dim_\kk 
\left(\fac{\kk[\ku{x}]\,}{I}\right)_q =
\left\{\begin{array}{ll}
1&\mbox{ if } q\in \CA(\pZ^n)\\
0& \mbox{ otherwise.}
\end{array}\right.
\]
Examples are the so-called toric ideal 
$J_\CA:=\big(x^a-x^b\,\big|\; a,b\in \pZ^n, \,a-b\in \ker\CA\big)$
and all its Gr\"obner degenerations. Indeed, these ideals form an irreducible,
the ``coherent'' component in the parameter space of all $\CA$-graded ideals.\\
The importance of $\CA$-graded ideals seems to be two-fold. First, they give
insight into a small layer of the deformation space of monomial ideals. Second,
via taking radicals and using the Stanley-Reisner construction, the monomial
$\CA$-graded ideals provide triangulations of the convex cone
$\CA(\pR^n)$. Hence, the set of those triangulations may be studied
by algebraic tools, cf.\ \cite{Diane}.
\par


\neu{intro-chains}
Sturmfels has investigated intrinsic properties of monomial
ideals that are satisfied for coherent ideals, but fail for general
$\CA$-graded ones. First, there is a combinatorial property
obtained by observing the
vertices of the fibers of $\CA$ restricted to $\pZ^n$. The second property
involves more algebraic concepts such as the degree of the generators
of the ideal, cf.\ \cite{agraded}, \cite{polytopes}.\\
In this context, Ho\c{s}ten and Thomas have 
addressed another point. In \cite{primechain} they observe
that monomial degenerations of toric ideals admit a very special primary
decomposition; the associated prime ideals occur in chains: 
\par

{\bf Definition:}
The ideal $I$ fulfils the {\em chain property} for its associated primes
if, for any associated, non-minimal prime $\kP$ of $I$, there is another one
$\kP^\prime\subseteq\kP$ with $\height(\kP^\prime)=\height(\kP)-1$.
\par

The subject of the present paper is to show that this property is of the 
same type as the first two mentioned above, 
i.e., not true for non-coherent $\CA$-graded monomial ideals, in general. 
In detail, we will prove the following
\par

{\bf Theorem \zitat{dim2}{dim2}/\zitat{dim3}{dim3}:}
{\em
The chain property holds for $\CA$-graded monomial ideals of dimension 
$d\leq 2$.
However, there are counter examples for $d=3$.  
}
\par


\neu{intro-primary}
The main tool for proving the previous theorem is the explicit knowledge of the
primary decomposition of monomial ideals from \cite{primary}.\\
However, instead of
quoting the result, we start our paper in \S \ref{primary}
with the presentation of a  
generalization to a certain class of binomial ideals, 
cf.\ Theorem \zitat{primary}{decomposition}. The fact that
primary decomposition does not leave the category of binomial ideals at all
follows already from \cite{binomial}.
\par

%
%
\section{Primary decomposition of saturated binomial ideals}\label{primary}

 
\neu{primary-T}
Let $I\subseteq \kk[\ku{x}]:=\kk[x_1,\dots,x_n]$ be a binomial ideal.
By $T(I)\subseteq\pZ^n$ we denote the set
\[
T(I):= \big\{a\in\pZ^n\,\big|\; x^a\notin I\big\}
\]
of the non-monomials in $I$. The set $T:=T(I)$ has the property
\kitem
\item[(i)]
If $a,b\in\pZ^n$ with $a\geq b$ (i.e., $a-b\in \pZ^n$) and
$a\in T$, then $b\in T$.
\kenditem
Every $T$ fulfilling this property occurs as $T(I)$ for some
binomial (even monomial) ideal. For the upcoming definition, we also need 
the following notation. If $\ell\subseteq [n]:=\{1,\dots,n\}$ is an arbitrary
subset, then we write
\vspace{-1ex}
\[
\Z^\ell:=\{a\in\Z^n\,|\; \supp\, a \subseteq \ell\}
\hspace{2em}\mbox{and}\hspace{2em}
\pZ^\ell:=\pZ^n\cap\Z^\ell\,.
\vspace{-3ex}
\]
\par


\neu{primary-standard}
{\bf Definition:}
Let $T\subseteq\pZ^n$ with property (i). A set $(r,\ell):=r+\pZ^\ell$
with $r\in\pZ^n$ and $\ell\subseteq [n]$ disjoint to 
$\supp\, r$ is called a {\em standard pair} if it is maximal (with respect
to inclusion) for the property $r+\pZ^\ell\subseteq T$. 
\par

Denoting by $B:=B(T)$ the set of standard pairs, we obviously have 
$T=\cup_{(r,\ell)\in B} (r+\pZ^\ell)$. 
The following result is well known, however, we include a short proof here
for the reader's convenience.
\par

{\bf Proposition:}
{\em
The set $B(T)$ is finite.
}
\par

{\bf Proof:}
Otherwise,
let $\ell\subseteq [n]$ be a subset such that there are infinitely many $r^i$
with $(r^i,\ell)\in B$. Then, we may choose a subsequence of $(r^i)$ 
that is increasing
via the partial order provided by $\pZ^{[n]\setminus\ell}\subseteq\pZ^n$:
If there is one value occurring infinitely often
as the entry of the $r^i$ in one of the $[n]\setminus\ell$ coordinates, then
this follows via induction by $\#([n]\setminus\ell)$. If not,
then the claim is trivial, anyway.
\vspace{1ex}\\
Assuming that the $r^i$ are increasing, then there is at least
one coordinate (e.g.\ the first one) that becomes arbitrary large.
Since $\big(r^1+(r^i_1-r^1_1)e^1\big)\leq r^i$, the property (i) implies
$\big(r^1+(r^i_1-r^1_1)e^1\big)+\pZ^\ell\subseteq T$; hence
$r^1+\pZ^{\ell\cup\{1\}}\subseteq T$. In particular, the set $r^1+\pZ^\ell$
was not minimal, i.e., $(r^1,\ell)\notin B$.
\qed
\par


\neu{primary-Tl}
Let $T\subseteq\pZ^n$ with property (i) of \zitat{primary}{T}. Then we define
for any subset $\ell\subseteq[n]$
\[
T(\ell) := \bigcup_{(r^i,\ell)\in B} \big(r^i+\pZ^\ell\big)\subseteq T
\]
and its closure via the partial order ``$\geq$''
\[
\ko{T(\ell)} :=  \bigcup_{(r^i,\ell)\in B,\, r\leq r^i}
\hspace{-0.5em}\big(r+\pZ^\ell\big)\subseteq T\,.
\vspace{-2ex}
\]
\par

{\bf Remark:}
$\;T(\ell)=\{r\in\pZ^n\,|\; r+\pZ^\ell\subseteq T\;
\mbox{and $\ell$ is maximal with this property}\}$.
\par


\neu{primary-saturated}
{\bf Definition:}
Let $I\subseteq \kk[\ku{x}]$ be a binomial ideal.
Define $K(I)\subseteq \Z^n$ as the abelian subgroup generated as
\[
K(I):=\big\langle a-b\,\big|\; 
a,b\in T(I) \;\mbox{with}\;
\lambda_a x^a - \lambda_b x^b \in I \;\mbox{for some}\, 
\lambda_a,\lambda_b\neq 0\big\rangle.
\]
The binomial ideal $I$ is called {\em saturated} if for any
$a,b\in T(I)$ with $a-b\in K(I)$ there are coefficients
$\lambda_a,\lambda_b\neq 0$ such that 
$\lambda_a x^a - \lambda_b x^b \in I$.
\par

{\bf Remark:}
The abelian subgroup $K(I)\subseteq \Z^n$ is minimal for $I$ to be
$\fac{\Z^n}{K(I)}$-homogeneous. Moreover, using this language, the 
saturation property means that the associated Hilbert function of the graded
ring $\fac{\kk[\ku{x}]\,}{I}$ always yields $0$ or $1$.
\par

{\bf Examples:}
\kitem
\item[1)] Monomial ideals: Here is $K(I)=0$.
\item[2)] $\CA$-graded ideals with $\CA:\Z^n\to\Z^d$:
We have $K(I)=\big\langle (a-b)\in \ker \CA\,|\; a,b\in T(I)\big\rangle
\subseteq \ker \CA$.\\
In particular, the map
$\fac{\Z^n}{K(I)}\surj \fac{\Z^n}{\ker\CA} \hookrightarrow\Z^d$ shows that the
$\Z^d$-grading is in general weaker than the $\fac{\Z^n}{K(I)}$-grading.
\item[3)] The ideal $I:=(x^2-xy)$ is {\em not} saturated: 
While $K(I)=(1,-1)\cdot\Z\subseteq\Z^2$, we have $x-y\notin I$.\hspace{-1em}
\kenditem
\par

{\bf Proposition:}
{\em
Let $I\subseteq \kk[\ku{x}]$ be a binomial ideal. 
If $I$ is saturated, then $T=T(I)$ 
and the fibers $T_q:=\{a\in T\,|\;a\mapsto q\}$ 
with $q\in \fac{\Z^n}{K(I)}$
fulfil, in addition to 
\zitat{primary}{T}{\rm(}i{\rm)}, the following property: 
\vspace{-0.5ex}
\kitem
\item[{\rm(ii)}]
If $g\in\pZ^n$ and $q\in\fac{\Z^n}{K(I)}$, then
$(g+T_q)\cap T=\emptyset$ or $g+T_q\subseteq T$.
\kenditem
}
\par

{\bf Proof:}
If $a,b\in T_q$, then there is an element 
$\lambda_a x^a - \lambda_b x^b \in I$, hence, 
$\lambda_a x^{a+g} - \lambda_b x^{b+g} \in I$.
Since $\lambda_a, \lambda_b\neq 0$, it follows that the latter monomials are
either both contained in $I$ or both not.
\qed
\par 

{\bf Open problem:}
Do there exist $\CA$-graded ideals $I$ containing at least one monomial, but
still satisfying $K(I)=\ker \CA$ or at least $K(I)_\Q=\ker\CA_\Q$?
\par

 
\neu{primary-TlLemma}
{\bf Lemma:}
{\em
Assume that $T\subseteq\pZ^n$ satisfies, for some subgroup 
$K\subseteq \Z^n$, the properties
{\rm(}i{\rm)} and {\rm(}ii{\rm)} of 
\zitat{primary}{T} and \zitat{primary}{saturated}, respectively.
\kitem
\item[\rm(a)]
If $a,b\in T_p$ and $g,h\in T_q$, then 
$a+g\in T$ iff $b+h\in T$.
\item[\rm(b)]
If $T(\ell)_q\neq\emptyset$ , then $T(\ell)_q=T_q$.
\item[\rm(c)]
If $\ko{T(\ell)}_q\neq\emptyset$ , then $\ko{T(\ell)}_q=T_q$.
\vspace{-1ex}
\kenditem 
In particular, both sets $T(\ell)\subseteq T$ and 
$\ko{T(\ell)} \subseteq T$ are unions
of selected whole fibers $T_q$.
}
\par

{\bf Proof:}
Part (a) uses property (ii) twice to compare 
$(a+g)$, $(a+h)$, and $(b+h)$ successively.
\vspace{0.5ex}\\
For (b) let $a\in T(\ell)_q$, $b\in T_q$. We will use 
Remark \zitat{primary}{Tl} several times: First, it follows that
$a+\pZ^\ell\subseteq T$. Then, with $g$ browsing through $\pZ^\ell$, property
(ii) implies that also $b+\pZ^\ell\subseteq T$. 
Moreover, the same argument applied in the opposite direction
shows that $\ell$ is maximal with this property, i.e., $b\in T(\ell)_q$.
\vspace{0.5ex}\\
Finally, let $a\in \ko{T(\ell)}_q$, $b\in T_q$. In particular, there is an
element
$g\in\pZ^n$ such that $(a+g)\in T(\ell)_{q+p}$ with $p$ being the image of $g$
via $\,\Z^n\to\fac{\Z^n}{K}$.
Applying the parts (a) and (b) successively, this means that 
$(b+g)\in T_{q+p}=T(\ell)_{q+p}$, hence $b\in \ko{T(\ell)}_q$.  
\qed
\par


\neu{primary-decomposition}
Let $I\subseteq \kk[\ku{x}]$ be a binomial ideal such that $T=T(I)$ and
$K=K(I)$ meet the conditions
{\rm(}i{\rm)} and {\rm(}ii{\rm)} of
\zitat{primary}{T} and \zitat{primary}{saturated}, respectively.
Then, for any $\ell\subseteq[n]$ such that some
$(\bullet,\ell)$ occurs as a standard pair, we define the ideal
\[
I^{(\ell)}\,:=\; I + \big(x^a\,|\; a\notin \ko{T(\ell)}\big) 
\;\subseteq\; \kk[\ku{x}]\,.
\]
If there was no standard pair containing $\ell$, then
$T(\ell)=\ko{T(\ell)}=\emptyset$, and the above definition would yield
$I^{(\ell)}=(1)$, anyway. 
\par

{\bf Theorem:}
{\em
If $I\subseteq \kk[\ku{x}]$ is a binomial ideal fulfilling 
{\rm(}i{\rm)} and {\rm(}ii{\rm)}, e.g.\ if $I$ is a saturated 
binomial ideal in the sense of \zitat{primary}{saturated},
then $\;I=\bigcap_{\ell} I^{(\ell)}$ will be a primary decomposition. 
}
\par

{\bf Proof:}
{\em Step 1: Being primary may be checked by means of the homogeneous
elements only:}\\
If $J\subseteq R$ is an ideal, then $J$ is primary if and only if the 
multiplication maps $(\cdot r):R/J\to R/J$ are either injective or nilpotent.
On the other hand, if $J$ is homogeneous in a graded ring, then these two
properties of linear maps $\psi:R/J\to R/J$ may be checked by using
homogeneous arguments only. Moreover, the sum of injective, homogeneous maps
of different degrees remains injective, the sum of nilpotent maps remains
nilpotent, and the sum of an injective and a nilpotent map is 
\vspace{1ex}
injective.\\
{\em Step 2: The ideals $I^{(\ell)}$ are primary:}\\
For $q\in\fac{\Z^n}{K}$ denote by 
$F_q:=\{a\in\pZ^n\,|\; a\mapsto q\}$ the whole fiber of
$q$; in particular, $T_q\subseteq F_q$.
If $I^{(\ell)}$ was not primary, then there would be elements 
$s\in \kk[\ku{x}]_p$ and $t\in \kk[\ku{x}]_q$ such that
$st\in I^{(\ell)}$, $s\notin I^{(\ell)}$, and $t^N\notin I^{(\ell)}$ for
every $N\geq 1$. Moreover, if we replace $s,t$ by
different representatives of their equivalent classes in
$\kk[\ku{x}]/I^{(\ell)}$, then the previous property does not
change.\\
For any degree $q\in\fac{\Z^n}{K}$ we know that
$I_q\subseteq I^{(\ell)}_q\subseteq \kk[\ku{x}]_q$ and 
$\,\dim_\kk \kk[\ku{x}]_q/I_q\leq 1$. 
Applied to our special situation this means that 
$I_p=I^{(\ell)}_p\subset \kk[\ku{x}]_p$, 
$I_{Nq}=I^{(\ell)}_{Nq}\subset \kk[\ku{x}]_{Nq}$, and $s,t$ may be assumed to
be monomials $x^a$ and $x^b$, respectively. Moreover, the product
$st=x^{a+b}$ is either contained in $I_{p+q}$, or we have that
$I^{(\ell)}_{p+q}= \kk[\ku{x}]_{p+q}$. Translated into the language of 
exponents, this means:
\[
a\in T_p=\ko{T(\ell)}_p
\hspace{0.6em}\mbox{and}\hspace{0.6em}
N b\in T_{Nq}=\ko{T(\ell)}_{Nq} \;\forall N\geq 1\,,
\hspace{0.6em}\mbox{but}\hspace{0.6em}
a+b\notin T
\hspace{0.6em}\mbox{or}\hspace{0.6em}
T_{p+q}\neq\ko{T(\ell)}_{p+q}\,.
\]
Since $\ko{T(\ell)}$ consists of only finitely many $\pZ^{\ell}$-slices, the
fact that $N b\in \ko{T(\ell)}$ for all $N\geq 1$ implies that
$b\in\pZ^{\ell}$. Hence, the property $a\in\ko{T(\ell)}$ yields 
$a+b\in \ko{T(\ell)}_{p+q}\subseteq T_{p+q}$ immediately.
Moreover, by Lemma \zitat{primary}{TlLemma}(c), the latter two sets have to
be equal, and we obtain a contradiction.
\vspace{1ex}\\
{\em Step 3: The intersection yields $I$:}\\
For every $q$ we have to show that there is at least one $\ell$ such that
$I^{(\ell)}_q=I_q$, i.e.\ such that $\ko{T(\ell)}_q=T_q$. However, 
if $T_q\neq\emptyset$,
the latter equality is equivalent to
$\ko{T(\ell)}\cap T_q\neq \emptyset$ 
by Lemma \zitat{primary}{TlLemma}(c). 
Hence, everything follows from 
$T=\bigcup_{\ell}T(\ell)=\bigcup_{\ell}\ko{T(\ell)}$.
\qed
\par

%
%
\section{Two-dimensional $\CA$-graded monomial ideals}\label{dim2}


\neu{dim2-monomial}
If $I\subseteq\kk[\ku{x}]$ is a {\em monomial} ideal, then 
Theorem \zitat{primary}{decomposition} yields the well known formula
for the primary decomposition of $I$ into the easier looking
$I^{(\ell)}=\big(x^a\,|\; a\notin \ko{T(\ell)}\big)$.
In particular, the associated primes are 
\vspace{-1ex}
\[
\kP^{(\ell)}=\sqrt{I^{(\ell)}}=\big(x_i\,\big|\; i\notin\ell\big)
\hspace{0.5em}\mbox{with $\ell$ such that there is a } (\bullet,\ell)\in B(T).
\]
%
If $I\subseteq\kk[\ku{x}]$ is additionally {\em $\CA$-graded} with respect to
some linear map $\CA:\Z^n\to\Z^d$ with $(\ker\CA)\cap\pZ^n=0$,
cf.~\zitat{intro}{agraded} for a definition,
then we know that $\CA$ induces an isomorphism
$T\stackrel{\sim}{\longrightarrow}\CA(\pZ^n)$.
In particular, every $\ell$ occurring in $B(T)$ fulfils
$\#\ell\leq d$, and the chain property, cf.~\zitat{intro}{chains}, looks as
follows:\\
$I$ fulfils the {\em chain property} for its associated primes if
for any non-maximal $\ell$ in $B(T)$ there is another 
$\ell^\prime\supseteq\ell$ in $B(T)$ with $\#\ell^\prime=\#\ell+1$.
\par

{\bf Remark:}
An $\ell$ occurring in $B(T)$ via some $(\bullet,\ell)$
is maximal if and only if $(0,\ell)\in B(T)$.
\par

We will show that the above chain property is always fulfilled 
for monomial, $\CA$-graded ideals as
long as $d\leq 2$, but it fails for $d\geq 3$.
\par


\neu{dim2-cap}
First, we need the following lemma describing the intersection behavior of two
different layers $(r,\ell)=(r+\pZ^\ell)$, $(s,m)=(s+\pZ^m)$ of the base
$B(T)$.
\par

{\bf Lemma:}
{\em
Let $T\subseteq\pZ^n$ be an arbitrary subset
satisfying the assumption {\rm(}i{\rm)} of \zitat{primary}{T}.
Then, two different $(r,\ell), (s,m)\in B(T)$ are either disjoint, or the
intersection is of the form
\[
(r+\pZ^\ell)\cap (s+\pZ^m) = (p+\pZ^{\ell\cap m})
\]
with strict inclusions $\ell\cap m\subset \ell,m$.
}
\par

{\bf Proof:}
Assume that the 
intersection is not empty. Then, outside $(\ell\cup m)$, 
the values of $r$ and $s$ coincide,
and we set $p$ as the common one. Within $(\ell\cup m)$ we define
\[
p_{|(\ell\setminus m)}:= s_{|(\ell\setminus m)}
\hspace{1em}\mbox{and}\hspace{1em}
p_{|(m\setminus \ell)}:= r_{|(m\setminus \ell)}.
\]
It remains to check that neither of $\ell, m$ is a subset of the other one. 
But if this was the case, say $\ell\subseteq m$, then 
$(r+\pZ^\ell)\subseteq (s+\pZ^m)$ implying that $(r+\pZ^\ell)$ would not
be a maximal subset of $T$, i.e., $(r,\ell)\notin B(T)$.
\qed
\par


\neu{dim2-triang}
{\bf Lemma:}
{\em
Let $I\subseteq\kk[\ku{x}]$ be a monomial, $\CA$-graded ideal. 
\kitem
\item[\rm(a)]
The set of all $\CA(\pR^\ell)$ with $\pZ^\ell\subseteq T$
forms a triangulation of the convex cone $\CA(\pR^n)\subseteq\R^d$.
The maximal cells come from those $\ell$ with $(0,\ell)\in B(T)$.
\item[\rm(b)]
Let $(0,\ell)\in B(T)$. Then, the map $\CA$ induces a natural bijection
\[
\{r\,|\; (r,\ell)\in B(T)\} 
\;\stackrel{\sim}{\longrightarrow}\;
\fac{\CA(\Z^n)\,}{\CA(\Z^\ell)}\;.
\]
\vspace{-3ex}
\kenditem
}
\par

{\bf Proof:}
Using Lemma \zitat{dim2}{cap} for part (a),
this and the injectivity in (b) are both 
simple consequences from the fact that the map $\CA$ is injective on the
subset $T\subseteq\pZ^n$.\\
To show the surjectivity in (b) we remark first that 
$\CA(\Q^\ell)=\CA(\Q^n)$. In particular, if some class
$\ko{w}\in\fac{\CA(\Z^n)\,}{\CA(\Z^\ell)}$ is represented by a
$w=\CA(a-b)$ with $a,b\in\pZ^n$, then there is an $N\geq 1$ with
$N\cdot\CA(b)\in\CA(\Z^\ell)$.
Hence $\ko{w}$ may be represented by 
$w+N\cdot\CA(b)=\CA\big(a+(N-1)b\big)$, 
i.e., by an element from $\CA(\pZ^n)$.\\
Let $\ko{w}$ be now represented from $\CA(\pZ^n)$. 
Since $\CA(\pZ^n)=\CA(T)$, this means that
\[
\forall\, a\in\pZ^\ell \hspace{0.8em} 
\exists\, \big(r(a),\ell(a)\big)\in B(T):
\hspace{0.8em} w+\CA(a)  \,\in\,
\CA\big(r(a)\big)+\CA\big(\pZ^{\ell(a)}\big).
\]
It remains to show that $\ell$ itself appears in the previous list of 
the sets
$\ell(a)$. But, if not, then all elements of $w+\CA(\pZ^\ell)$ would be
contained in at most finitely many shifts of maximal cells different from 
$\CA(\pZ^\ell)$. 
\qed
\par


\neu{dim2-dim2}
{\bf Theorem:}
{\em
Let $I\subseteq\kk[\ku{x}]$ be a monomial, $\CA$-graded ideal of dimension
$d\leq 2$. Then $I$ satisfies the chain property for its associated primes.  
}
\par

{\bf Proof:}
Since there is nothing to show for one-dimensional ideals, we consider the 
case of $d=2$. Let us assume that the chain condition is violated, i.e., 
for every $(r,\ell)\in B(T)$ we have either
$\#\ell=2$ or $\ell=\emptyset$, and, moreover, there is at least one
$(r^\ast,\emptyset)\in B(T)$ of the second type.\\
Since the $\ell$'s with cardinality two provide a triangulation of the
two-dimensional cone $\CA(\pR^n)\subseteq\R^2$, we may order them in a natural
way as $\ell^1,\dots,\ell^N$. Then Lemma \zitat{dim2}{cap} implies that
adjacent sets $\ell^{i-1}, \ell^i$ share a common element, say $i$. 
Denoting the canonical basis elements of $\Z^n$ by $e^i$, this yields
the following setup
\[
\ell^i=\{i,i+1\}
\hspace{0.7em}\mbox{and}\hspace{0.7em}
\sigma_i:=\CA\big(\pR^{\ell^i}\big)=
\big\langle \CA(e^i), \CA(e^{i+1})\big\rangle_{\pR}
\hspace{0.5em}\mbox{with}\hspace{0.7em} i=1,\dots,N (<n).
\]
\begin{center}
{
\unitlength=0.6pt
\begin{picture}(444.00,301.00)(0.00,0.00)
\put(349.00,17.00){\makebox(0.00,0.00){$f^1$}}
\put(314.00,100.00){\makebox(0.00,0.00){$f^2$}}
\put(258.00,171.00){\makebox(0.00,0.00){$f^3$}}
\put(103.00,260.00){\makebox(0.00,0.00){$f^N$}}
\put(-3.00,253.00){\makebox(0.00,0.00){$f^{N+1}$}}
\put(9.00,271.00){\vector(-4,-1){31.00}}
\put(121.00,269.00){\vector(-4,1){32.00}}
\put(286.00,168.00){\vector(-1,1){19.00}}
\put(341.00,89.00){\vector(-1,3){10.00}}
\put(366.00,1.00){\vector(0,1){23.00}}
\put(10.00,302.00){\makebox(0.00,0.00){$\CA(e^{N+1})$}}
\put(123.00,301.00){\makebox(0.00,0.00){$\CA(e^N)$}}
\put(341.00,215.00){\makebox(0.00,0.00){$\CA(e^3)$}}
\put(417.00,108.00){\makebox(0.00,0.00){$\CA(e^2)$}}
\put(450.00,1.00){\makebox(0.00,0.00){$\CA(e^1)$}}
\put(140.00,117.00){\makebox(0.00,0.00){$\ddots$}}
\put(71.00,173.00){\makebox(0.00,0.00){$\sigma_N$}}
\put(213.00,74.00){\makebox(0.00,0.00){$\sigma_2$}}
\put(247.00,25.00){\makebox(0.00,0.00){$\sigma_1$}}
\put(177.00,160.00){\makebox(0.00,0.00){}}
\put(166.00,146.00){\makebox(0.00,0.00){}}
\put(168.00,153.00){\makebox(0.00,0.00){}}
\put(77.00,1.00){\line(1,6){48.00}}
\put(77.00,1.00){\line(-1,4){72.00}}
\put(77.00,1.00){\line(5,4){252.00}}
\put(77.00,1.00){\line(3,1){306.00}}
\put(77.00,1.00){\line(1,0){345.00}}
\put(520.00,50.00){\makebox(0.00,0.00){Figure 1}}
\end{picture}}
\end{center}
Using part (b) of Lemma \zitat{dim2}{triang}, 
we may choose, for every $i$, a pair
$(r^i,\ell^i)\in B(T)$ such that
$A(r^\ast)\in\CA(r^i)+\CA(\Z^{\ell^i})$.
On the other hand, $r^\ast$ is ``isolated'', i.e., it does not belong to any of
the sets $r^i+\pZ^{\ell^i}$. Then the injectivity of $\CA$ on $T$ implies that
$\CA(r^\ast)\notin \CA(r^i)+\CA\big(\pZ^{\ell^i}\big)$, and since 
$\CA\big(\pZ^{\ell^i}\big) = 
\CA\big(\Z^{\ell^i}\big) \cap \CA\big(\pR^{\ell^i}\big)$,
we obtain $\CA(r^\ast)\notin \CA(r^i)+\sigma_i$.
\vspace{0.3ex}
If $f^i$ denote linear forms on $\Z^2$ such that the cones $\sigma_i$ 
are given by
\vspace{-1ex}
\[
\sigma_i=\big\{x\in\R^2\,\big|\; \langle x, f^i\rangle \geq 0,\;
\langle x, f^{i+1}\rangle \leq 0\big\}\,,
\]
cf.\ Figure 1, then this statement can be translated into
\[
\big\langle \CA(r^\ast) - \CA(r^i),\, f^i\big\rangle <0
\hspace{1.2em}\mbox{or}\hspace{1.0em}
\big\langle \CA(r^\ast) - \CA(r^i),\, f^{i+1}\big\rangle >0
\hspace{1em}\mbox{for }\,i=1,\dots,N.
\]
Reorganizing these inequalities, we obtain that either
\[
\big\langle \CA(r^\ast),\, f^1\big\rangle <
\big\langle \CA(r^1),\, f^1\big\rangle,
\hspace{0.9em}\mbox{or}\hspace{0.9em}
\big\langle \CA(r^N),\, f^{N+1}\big\rangle <
\big\langle \CA(r^\ast),\, f^{N+1}\big\rangle,
\]
or there is an $i\in\{2,\dots,N\}$ with 
\[
\big\langle \CA(r^{i-1}),\, f^i\big\rangle <
\big\langle \CA(r^\ast),\, f^i\big\rangle <
\big\langle \CA(r^i),\, f^i\big\rangle
\]
as depicted in Figure 2.
\begin{center}
{
\unitlength=0.6pt
\begin{picture}(427.00,348.00)(0.00,0.00)
\put(125.00,57.00){\makebox(0.00,0.00){$\CA(r^\ast)+\sigma_{i-1}$}}
\put(70.00,111.00){\makebox(0.00,0.00){$\CA(r^\ast)+\sigma_i$}}
\put(140.00,170.00){\circle*{8.00}}
\put(158.00,97.00){\circle*{8.00}}
\put(19.00,19.00){\circle*{8.00}}
\put(355.00,280.00){\makebox(0.00,0.00){$(f^i)^\bot$}}
\put(0.00,0.00){\makebox(0.00,0.00){$\CA(r^\ast)$}}
\put(19.00,19.00){\line(1,2){164.00}}
\put(19.00,19.00){\line(4,3){330.00}}
\put(19.00,19.00){\line(6,1){410.00}}
\put(325.00,157.00){\makebox(0.00,0.00){$\CA(r^{i-1})+\sigma_{i-1}$}}
\put(250.00,295.00){\makebox(0.00,0.00){$\CA(r^i)+\sigma_i$}}
\put(140.00,170.00){\line(1,2){110.00}}
\put(140.00,170.00){\line(4,3){230.00}}
\put(158.00,97.00){\line(4,3){230.00}}
\put(158.00,97.00){\line(6,1){270.00}}
\put(450.00,30.00){\makebox(0.00,0.00){Figure 2}}
\end{picture}}
\end{center}
Assume, w.l.o.g., that the latter two inequalities apply to some $i$.
Then we consider the series 
$\CA(r^\ast)+\Z_{\geq 0}\CA(e^i)\subseteq \CA(r^\ast)+(f^i)^\bot$:
All its 
members are contained in $\CA(T)$. Hence, up to finitely many exceptions,
they have to be contained in some $\CA(s)+\CA(\pZ^\ell)$ with 
$(s,\ell)\in B(T)$ and $\#\ell=2$.\\
On the other hand, if 
$\ell \neq \ell^{i-1},\ell^i$, then any shift of the cone
$\CA\big(\pR^\ell\big)$ intersects the ray $\CA(r^\ast)+\R_{\geq 0}\CA(e^i)$
in a compact set, i.e., the intersection contains 
at most finitely many lattice points. 
Thus, almost all elements of $\CA(r^\ast)+\Z_{\geq 0}\CA(e^i)$ are contained in
sets of the form 
$\CA(s)+\CA(\pZ^{\ell^{i-1}})$ or $\CA(s)+\CA(\pZ^{\ell^i})$ with
$(s,\ell^{i-1}), (s,\ell^i)\in B(T)$.
\vspace{0.5ex}\\
However, applying part (b) of Lemma \zitat{dim2}{triang} again, we see that 
there is no freedom left for the element $s$. 
It has to equal $r^{i-1}$ or $r^i$, respectively. But since the sets
$\CA(r^{i-1})+\sigma_{i-1}$ and $\CA(r^i)+\sigma_i$ do not meet the ray
$\CA(r^\ast)+\R_{\geq 0}\CA(e^i)$ at all, we have obtained a contradiction.
\qed
\par

%
%
\section{A counter example in dimension three}\label{dim3}


\neu{dim3-dim3}
Roughly speaking, the proof of the previous theorem 
worked as follows: We have shown that
the shifted two-dimensional cells of the triangulation create gaps that cannot
be filled with isolated $T$-elements only. In dimension three this concept
fails, since some cells might be arranged in cycles. In particular, we have
\par

{\bf Theorem:}
{\em
There exists an example of a monomial, $\CA$-graded ideal
$I\subseteq\kk[\ku{x}]$ with $d=3$ such that the chain property for the
associated primes is violated.
}  
\par

{\bf Proof:}
Take $n=16$ with the variables
$e_i,f_i,g_i$ ($i=1,2,3$) and $k_\nu$ ($\nu=1,\dots,7$). 
The ideal $I$ is defined by the following 100 generators
\kitem
\item
$f_i\,k_\nu,\;
g_i\,k_\nu,\;
g_i\,g_j,\;
k_\nu\, k_\mu$
with $i,j\in\{1,2,3\}$ and $\mu,\nu\in\{1,\dots,7\}$,
\item
$f_i^2,\;
f_i\,g_{i+1},\;
f_{i-1}\,f_i\,g_{i-1}$,
\item
$f_i\,g_i\,e_{i-1},\;
f_{i+1}\,g_i\,e_{i-1},\;
f_i\,f_{i+1}\,e_{i-1}$, and
\item
$f_i\,e_{i-1}\,e_{i+1},\;
g_i\,e_{i-1}\,e_{i+1}$
\kenditem
with $i\in\fac{\Z}{3\Z}$. Denoting by $\{E_1,E_2,E_3\}$ the canonical basis 
$\{(1,0,0), (0,1,0), (0,0,1)\}$ of
$\Z^3$, then $I$ is $\CA$-graded with respect to 
the linear map $\CA:\Z^{16}\to\Z^3$ given by
\[
e_i\mapsto 2\,E_i,\hspace{1em}
f_i\mapsto E_i,\hspace{1em}
g_i\mapsto E_i + (1,1,1),
\vspace{-3ex}
\]
and
\vspace{-1ex}
\[
\big\{k_1,\dots,k_7\big\}
\hspace{0.5em}\mapsto \hspace{0.5em}
\big\{ (3,2,2),\, (2,3,2),\, (2,2,3),\, (2,3,3),\, (3,2,3),\, 
(3,3,2),\, (3,3,3)\big\}\,.
\]
\vspace{-3ex}
\qed
\par

{\bf Remark:}
Since $\CA:\pZ^{16}\to\pZ^3$ is surjective, the associated toric ideal defines
the semigroup algebra $\kk\big[\pZ^3\big]$, 
i.e., the associated toric variety is $\kk^3$.
\par


\neu{dim3-example}
In addition to the plain presentation of the example in the previous section, 
we would also
like to show how it really works. In particular, we rather describe the set
$T=T(I)$ of standard monomials via its basis $B(T)$ -- and how these sets map
to $\Z^3$.\\ 
To improve the readability, we write $\langle\ell\rangle$ 
for the semigroup
$\pZ^\ell$. Then the ideal $I$ has the following maximal (with respect to the
partial order ``$\geq$'' induced by $\pZ^{16}$) standard pairs
\kitem
\item[(i)]
$k_\nu + \big\langle e_1,e_2,e_3\big\rangle$ ($\nu=1,\dots,7$),
\item[(ii)]
$(f_i+f_{i+1}) + \big\langle e_i, e_{i+1}\big\rangle, \hspace{0.5em}
(f_i+g_i) + \big\langle e_i, e_{i+1}\big\rangle, \hspace{0.5em}
(f_{i+1}+g_i) + \big\langle e_i, e_{i+1}\big\rangle$,
\item[(iii)]
$g_{i+1} +  \big\langle e_i, e_{i+1}\big\rangle$, \hspace{0.3em}and
\item[(iv)]
$(f_1+f_2+f_3)$
\kenditem
with $i\in\fac{\Z}{3\Z}$.
Dropping the maximality condition, we have to add the standard pairs
\kitem
\item[(v)]
$\big\langle e_1,e_2,e_3\big\rangle\;$ and
\item[(vi)]
$f_i +  \big\langle e_i, e_{i+1}\big\rangle, \hspace{0.5em}
f_{i+1} +  \big\langle e_i, e_{i+1}\big\rangle, \hspace{0.5em}
g_i +  \big\langle e_i, e_{i+1}\big\rangle$.
\kenditem
The violation of the chain property for the associated primes is caused by the
existence of the standard pair
$\big((f_1+f_2+f_3),\,\emptyset\big)$. The remaining standard pairs involve
only sets $\ell\subseteq\{1,2,\dots,16\}$ with $\#\ell\geq 2$.\\
Finally, for the $\CA$-graded property, let us consider the $\CA$-images:
\kitem
\item
(i) and (v) yield all triples with entries $0$ or $\geq 2$.
\item
Assuming $[i=1]$, the series (ii), (iii), (vi) provide
$\big(2\,\pZ, 2\,\pZ, 0\big)$ shifted by 
$(1,0,0)$, $(0,1,0)$, $(1,1,0)$, or $(3,1,1)$, $(2,1,1)$, $(1,2,1)$,
$(2,2,1)$. This means that every triple having $0$ or $1$ as its last entry is
reached, except for 
$\big(1, 1+2\,\pZ, 1\big)$ and $\big(2\,\pZ, 2\,\pZ, 0\big)$ itself.\\
However, the latter series already occurred in the previous point 
(i)/(v), and,
beginning with $(1,3,1)$, the first series is included in
(ii)/(iii)/(vi) with $[i=2]$.
\item
The isolated (iv) yields the missing triple $(1,1,1)$.
\kenditem
It follows that $\CA$ maps $T\subseteq\pZ^{16}$ onto $\pZ^3$. 
Moreover, a closer
look shows that the restriction $\CA_{|T}$ is injective, indeed.
\par

%
%

{\small
\parbox[t]{8cm}{
Klaus Altmann\\
Institut f\"ur Mathematik\\
Mathematisch Naturwissenschaftliche Fakult\"at II\\
Humboldt-Universit\"at zu Berlin\\
Rudower Chaussee 25\\
D-10099 Berlin, 
Germany\\
E-mail: {\tt altmann@mathematik.hu-berlin.de}}}

\end{document}